\documentclass[reqno,12pt]{amsart}
 
\newtheorem{theorem}{Theorem}[section]
\newtheorem{lemma}[theorem]{Lemma}
\newtheorem{corollary}[theorem]{Corollary}

\theoremstyle{definition}

\newtheorem{hypothesis}[theorem]{Hypothesis}

\theoremstyle{remark}
\newtheorem{remark}[theorem]{Remark}

\newcommand\bR{\mathbb{R}}

\newcommand\cC{\mathcal{C}}

\newcommand\cR{\mathcal{R}}

\newcommand\frC{\mathfrak{C}}

\newcommand{\mysection}[1]{\section{#1}
\setcounter{equation}{0}}

\newcommand\cbrk{\text{$]$\kern-.15em$]$}}
\newcommand\opar{
\text{\,\raise.2ex\hbox{${\scriptstyle |}$}\kern-.34em$($}}
\newcommand\cpar{%
\text{$)$\kern-.34em\raise.2ex\hbox{${\scriptstyle |}$}}\,}
\newcommand\obrk{\text{$[$\kern-.15em$[$}}

\begin{document}

\title[PDEs with growing coefficients]
{Elliptic and parabolic second-order PDEs with growing coefficients}

\author{N. V. Krylov}
\address{127 Vincent Hall \\ University of Minnesota,
Minneapolis, MN 55455  USA} \email{krylov@math.umn.edu}

\author{E. Priola}
\address{Dipartimento di Matematica\\Universit\`a di Torino
\\Via Carlo Alberto
10\\10123 Torino\\Italy } \email{ enrico.priola@unito.it}

\date{}

\thanks{The first author was partially supported by
NSF Grant DMS-0653121. The second named author gratefully
acknowledges the support by the  M.I.U.R. research projects Prin
2004 and 2006 ``Kolmogorov equations''.}



\begin{abstract} We consider a second-order parabolic equation in
$\bR^{d+1}$ with  possibly unbounded lower order coefficients. All
coefficients are
   assumed
to be  only  measurable   in the time variable
    and locally
H\"older continuous in the space variables.
   We show that   global  Schauder
   estimates  hold even in this case. The proof introduces
    a new localization procedure.
  Our  results  show that the constant appearing
  in the classical
 Schauder estimates is
  in fact   independent of the $L_{\infty}$-norms of the lower order
  coefficients.  
     We also give a proof of   uniqueness
      which is of independent interest
     even in the case  of bounded coefficients.

\end{abstract}

\maketitle

\mysection{Introduction}

   Let us consider the following second-order
operator $L$,
\begin{equation}                                    \label{lu}
L u = a^{ij}(t,x)u_{x^{i}x^{j}}(t,x) + b^{i}(t,x) u_{x^i}(t,x)
-c(t,x)u(t,x),
\end{equation}
   acting on functions defined
on  $[T, \infty) \times \bR^d$ if $T \in (-\infty,  \infty)$ and
   on $\bR^{d+1}$ if $T=-\infty$
(the summation convention is enforced throughout the article).
   We prove   global Schauder estimates for
solutions of the equation
\begin{equation}
                                                \label{2.6.4}
u_{t}(t,x) + Lu(t,x) = f(t,x),\quad (t,x) \in ( T,   \infty)
\times
   \bR^{d}.
\end{equation}

Roughly speaking, we will assume that $a,b,c,f$ are measurable in
$t$ locally bounded with respect to $(t,x)$, H\"older continuous
in $x$ in any ball of radius one with a constant independent of
$t$ and the position of the ball.   Moreover,  we assume that
$c(t,x)$ is always greater than some  constant $\delta
 > 0$, $f$ is pointwise ``controlled" by $c$, and the matrix $a$
is uniformly bounded and uniformly positive definite. The local
H\"older continuity does not prevent $b,c$, and $f$ from growing
linearly as $|x|\to\infty$. Thus, we do not assume that  $b$, $c$,
and $f$ are globally bounded as in the classical setting (see
\cite{B}, \cite{Kn}, \cite{LSU}, \cite{L}, \cite{Lo}).
  Recently, the interest in elliptic and
parabolic equations with unbounded coefficients in the whole space
has  increased (see, for instance, \cite{BL}, \cite{Ce},
\cite{DL},  \cite{DG}, \cite{Lu1},
\cite{LV},
 \cite{MPW},  \cite{P}
and the references therein). Such equations arise naturally also
in stochastic control and filtering theory (see, for instance,
\cite{FM} and \cite{S}).

   In Theorem \ref{theorem 1.10.1},  we
obtain parabolic   Schauder estimates of the type
\begin{equation}
                                                    \label{brandt}
\sup_{t\ge  T} \| u(t,\cdot) \|_{ 2+\alpha}\leq N \sup_{t\ge  T}
\| f(t,\cdot)\|_{\alpha,\, loc},
\end{equation}
by means of a new localization procedure   and Lemma 
\ref{lemma 1.10.1} saying that the constants in
classical Schauder estimates
for equations with coefficients depending only on $t$
are independent of the magnitudes of $b^i$ and $c$. In particular, from
\eqref{brandt} we deduce new elliptic Schauder estimates  when $a^{ij}$,
$b^i$, $c$,
   and $f$ do not
depend on $t$.     It is noteworthy that  to prove
    the new
   elliptic Schauder estimates in $\bR^d$ we need to
use the corresponding result for {\em parabolic\/} equations
   in $\bR^{d+1}$.
Estimate \eqref{brandt} allows us to prove the solvability of
\eqref{2.6.4} (Theorem \ref{theorem 2.7.1}), of the
   related  Cauchy problem
(Theorem \ref{theorem 4.15.3}), and of similar elliptic equations
in the whole space (Theorem \ref{5.7.2}).

While dealing with equations with growing coefficients in the
whole space it is quite natural to work in $C^{2+\alpha}$ spaces
with weights which would  not allow the derivatives of solutions
    to grow   so that the
terms $a^{ij}u_{x^{i}x^{j}}$, $b^{i}u_{x^{i}}$, and $cu$ would
remain in the usual $C^{\alpha}$ without weights (see, for
instance,
    \cite{CV},
\cite{Krk}, and the references therein).

   Naturally, once we want to allow
    the  coefficients to grow, the question arises as to what
happens if the coefficients do not grow but $f$ does. Such cases
were investigated for instance
   in \cite{Lo} (see also \cite[Remark 2.2]{Kry2})
where solutions
   were looked for in $C^{2+\alpha}$ spaces with weights.

   We discuss now some recent papers dealing with Schauder estimates
   for   elliptic equations and
    for {\it autonomous} parabolic Cauchy
problems involving  unbounded
    coefficients.

In \cite{CV}  elliptic and  parabolic equations with
   unbounded coefficients are studied
   assuming a kind of
   ``balance'' between the first order
   term $b^i u_{x^i}$ and the potential term $c u$ (a similar balance
   was
   also used
   in \cite{AB}).
    Schauder estimates in \cite{CV}  follow
    by
   generation of analytic semigroup in H\"older spaces.
   Weighted H\"older spaces $C^{2+\alpha}(V)$ are introduced such
that $L$ becomes a bounded operator from $C^{2+\alpha}(V)$ onto
the usual $C^{\alpha}$ with bounded inverse. Roughly speaking, $V$
in \cite{CV} is a function comparable with $c$ such that $|b|\leq
V^{1/2}$, and $C^{2+\alpha}(V)$ is the space of functions such
that $Vu$ and $u_{xx}$ belong to the usual $C^{\alpha}$.

    It was discovered in \cite{DL}  that
   even if   $b$  grows linearly  and $c=0$ one can still have the
   solvability theory in the {\em usual\/}   $C^{2+\alpha}$ spaces
without
   weights (note that in this case, contrarily to \cite{CV},
   generation of
    analytic semigroup fails).
   In \cite{DL} it is  assumed that $a^{ij}$ are constant  and
    $b(t,x) = A x$, $x \in \bR^d$, for some fixed
   $d \times d$ real  matrix $A$.

The results  of  \cite{DL}
    are surprising for the
following reason. For elliptic equations with bounded  H\"older
continuous coefficients we always have $LC^{2+\alpha}=C^{\alpha}$
and $C^{2+\alpha}=L^{-1}C^{ \alpha}$. However, if we allow $b$ to
grow, then $L^{-1}C^{ \alpha}\subset C^{2+\alpha}$ with proper
inclusion and, even more than that, for different $L$ the sets
$L^{-1}C^{ \alpha}$ may be different (we say more
about this in Remark \ref{new}).
   This situation differs dramatically from what was common before.
For instance, in \cite{CV} for all operators $L$, satisfying the
conditions imposed there with the same $V$, the set $L^{-1}C^{
\alpha}$ is always $C^{2+\alpha}(V)$.

   After \cite{DL}, several authors  dealt with elliptic Schauder
estimates  and Schauder estimates for autonomous   Cauchy problems
with unbounded {\it smooth\/} coefficients, see, for instance,
\cite{BL},  \cite{Ce0},   \cite{Ce}, \cite{Lu1}, \cite{LV},
\cite{P}.
   These papers also contain regularity results not covered
   in the present article. However,
our understanding is that their methods and results can not be
used to derive Schauder
   estimates in
   our situation (even if we would consider our coefficients
    $a$, $b$, and $c$ independent on $t$).
      For example, in \cite{Lu1} and \cite{Ce}
        Schauder estimates are proved assuming that $a^{ij}$ and $b^i$
          are smooth
           enough (since the methods used require to differentiate
             the
equations three times).  Note that
    mollifying the coefficients $a^{ij}$
   and then using , for instance, the results of \cite{Lu1}
   do  not
   seem to allow one to get  our elliptic Schauder estimates.
   In \cite{LV}   Schauder
   estimates are proved assuming only that  $a^{ij}$ are bounded
   and  have  first order
   bounded derivatives but  an additional  compatibility
     condition
   between  $a^{ij}$ and $b^i$ is imposed.

Classical  (possibly non-autonomous)
   parabolic Schauder
estimates when $a^{ij}$, $b^i$, and $c$ are bounded and H\"older
continuous in space and time are proved in \cite{LSU}. Partial
Schauder estimates like \eqref{brandt} in the case of bounded
coefficients which are only measurable in time were discovered in
\cite{B}. Then in \cite{Kn} it was   shown that the second
derivatives of the solution $u$ in $x$ are H\"older continuous
with respect to the time variable.  We will see below in Lemma
\ref{lemma 4.13.2} that this result is, actually, an embedding
theorem and has little to do with parabolic equations. In \cite{L}
interior parabolic estimates and Schauder estimates in a bounded
parabolic domain up to the boundary are proved.  The parabolic
Cauchy  problem when the coefficients are bounded, discontinuous
    in time and  H\"older continuous in space
is investigated in \cite{Lo}.

   It seems that our Theorem \ref{theorem 2.7.1}
    on solvability of equation \eqref{2.6.4}
    has not
   been stated
   before even in the case of bounded coefficients.
   We also stress
    that  uniqueness  of solutions is based
   on an a priori
   estimate of independent interest
   (see Theorem
\ref{theorem 4.25.1}) which seems to be new even
   when the coefficients are bounded.
    Our  exposition
   is independent
of the standard $C^{2+\alpha}$-theory of elliptic and parabolic
equations and, actually, one can get basic Schauder estimates of
this theory directly  from our results  (see Corollary
\ref{corollary 4.15.1}).

Our main results are stated in Section \ref{section 4.16.1} where
we also prove two of them, namely Theorem \ref{5.7.2} and
\ref{theorem 4.15.3}. Theorems \ref{theorem 1.10.1} and
\ref{theorem 2.7.1} are proved in Section \ref{section 4.16.3}
after we prepare necessary auxiliary results on equations with
coefficients independent of $x$ in Section \ref{section 4.16.2}.

\mysection{Main results}
                                           \label{section 4.16.1}

Introduce

$$
\bR^{d+1}=\{(t,x):t\in\bR,x=(x^{1},...,x^{d})\in\bR^{d}\}.
$$

\begin{hypothesis}                              \label{hy}

(i) The matrix $(a^{ij}(t,x))$ is symmetric and there exist
constants $\delta,K\in(0,\infty)$ such that
$$
KI \geq(a^{ij}(t,x))\geq\delta I , \quad c(t,x)\geq\delta,\quad
\;\;\; (t,x) \in\bR^{d+1} .
$$

(ii) The functions $a^{ij}(t,x)$, $b^i(t,x)$, and $c(t,x)$  are
measurable in $\bR^{d+1} $; $b(t,0)$ and $c(t,0)$ are  locally
bounded in $t$, and, for some $\alpha \in (0,1)$,
$$
|a^{ij}(t,x)-a^{ij}(t,y)| +
   |b^i(t,x)-b^i(t,y)
|+ |c(t,x)-c(t,y)|\leq K|x-y|^{\alpha}
$$
for all $t \in\bR $, $i,j=1,...,d$, and $x,y\in\bR^{d}$ such that
\begin{equation}
                                                  \label{4.19.1}
|x-y|\leq1.
\end{equation}
\end{hypothesis}

\begin{hypothesis}                              \label{H2}
The function $f(t,x)$
   is measurable in
   $\bR^{d+1} $   and there exist
    constants $F_{0}$ and $F_{\alpha}$
   such that we have
\begin{equation}
                                                  \label{5.7.3}
|f(t,x)|\leq F_{0}\, c(t,x)
\end{equation}
$$
|f(t,x)-f(t,y)|\leq F_{\alpha}|x-y|^{\alpha},
$$
whenever $t\in\bR$, $x,y\in\bR^{d}$  and \eqref{4.19.1} holds.
\end{hypothesis}

  Our Theorem \ref{theorem 4.15.3} about the Cauchy problem
 is  also true when $c(t,x)$ is
only nonnegative  on $\bR^{d+1}$. More precisely, such result
holds if $c(t,x)$ is replaced by  1 in  assumption \eqref{5.7.3}
with the other hypotheses  unchanged (see Remark \ref{remark 6.13.1}).

As usual $C^{\alpha}=C^{\alpha}(\bR^{d})$ is the Banach space of
real functions $g$ on $\bR^{d}$ with finite norm
$$
\|g\|_{\alpha}=\|g\|_{0}+[g]_{\alpha},
$$
where
$$
\|g\|_0=\sup_{x \in \bR^{d}}|g(x)|,\quad
[g]_{\alpha}=\sup_{x,y\in\bR^{d}}
\frac{|g(x)-g(y)|}{|x-y|^{\alpha}}\,,\quad\big(\mbox{\footnotesize{$
\frac{0}{0}:=0$}}\big).
$$
By $C^{2+\alpha}$ we denote the space of real functions $u$ on
$\bR^{d}$ with finite norm
$$
\|u\|_{2+\alpha}=\|u\|_{0}+\|Du\|_{0} +\|D^{2}u\|_{0} + [
u]_{2+\alpha},
$$
where $Du=(D_{1}u,...,D_{d}u)$, $D^{2}u=(D_{ij}u;i,j=1,...,d)$,
$$
D_{i}u=u_{x^{i}},\quad D_{ij}u=D_{i}D_{j}u=u_{x^{i}x^{j}},\quad
   [u]_{2+ \alpha} = [D^{2}u]_{\alpha} .
$$

   For any $T\in(-\infty,\infty)$,  we  define
$$
\bR^{d+1}_{T}= \bR^{d+1}\cap\{t\geq T\}= [T,  \infty) \times
\bR^d,
$$
and, if $T=-\infty$, we set $\bR^{d+1}_{T}=\bR^{d+1}$.
The supremum norm of functions on
$\bR^{d}_{T}$ will be denoted by $\|\cdot\|_{0,T}$.

We will be working with the set $\cC^{2+\alpha}(T)$ of functions
$u(t,x)$ defined on $\bR^{d+1}_{T}$ such that

(i) the function $u$ is continuous in $\bR^{d+1}_{T}$;

(ii) for each  finite $t\geq T$, we have $u(t,\cdot)\in
C^{2+\alpha} $ and $\|u(t,\cdot)\|_{2+\alpha}$ is bounded
in $t$;

(iii) there is a  measurable  function $g(t,x)$ defined on
$\bR^{d+1}_{T}$ such that for any $\zeta\in
C^{\infty}_{0}(\bR^{d+1})$ the function $\zeta(t,x) g(t,x)$ is
bounded  and $\alpha$-H\"older continuous in $x$ with constant
independent of $t$   and for any $x\in\bR^{d}$   and any finite
$s$ and $t$, such that  $T\leq s<t $, we have
\begin{equation} \label{g1}
u(t,x)-u(s,x)=\int_{s}^{t}g(r,x)\,dr.
\end{equation}

For such a function $u$ we denote $u_{t}=g$. Obviously $g$ is the
generalized derivative of $u$ with respect to $t$. By
$\cC^{2+\alpha}_{0}(T)$ we denote a subset of $\cC^{2+\alpha}(T)$
consisting of functions vanishing for large  $|x|+|t|$.

   Recall that a continuous function  $u$ has a bounded
generalized derivative with respect to a coordinate if and only if
it is Lipschitz continuous with respect to this coordinate and if
and only if $u$ is absolutely continuous with respect to the
coordinate and its classical derivative (existing almost
everywhere) is bounded. Under any of the above conditions the
classical derivative coincides with the generalized one and its
essential supremum equals the Lipschitz constant.

Accordingly,  solutions  of equation \eqref{2.6.4} will be looked
for in $\cC^{2+\alpha}(T)$ and in this class \eqref{2.6.4} is
equivalent to the fact that, for any $x \in \bR^d $  and finite
$t>s$ with $s\geq T$, we have
\begin{equation}
                                         \label{2.7.5}
u(t,x)-u(s,x)=\int_{s}^{t}f(r,x)\,dr- \int_{s}^{t}Lu(r,x)\,dr.
\end{equation}

\begin{remark}
Note that,  according to  Corollary \ref{corollary 4.15.3}, the
unique solution $u$ of \eqref{2.6.4} will have  more regularity.
\end{remark}

   Here are
   our  main results in which we always suppose
that Hypotheses \ref{hy} and \ref{H2} are satisfied  and
$T\in[-\infty,\infty)$.

The first result bears on an a priori estimate.

\begin{theorem}
                                       \label{theorem 1.10.1}

Let $u\in\cC^{2+\alpha} (T)$ satisfy \eqref{2.6.4}. Then:

(i) There is a constant $N=N(\delta,\alpha,d,K)$ such that, for
all finite $t\geq T$, we have
\begin{equation}
                                                       \label{4.15.3}
   \|u(t,\cdot)\|_{2+\alpha} \leq N (F_{0}+F_{\alpha}).
\end{equation}

(ii) If $f$ and the coefficients of $L$ are independent of $t$,
then with the same $N$, for any $u\in C^{2+\alpha} $,
$$
   \|u\|_{2+\alpha}\leq N   ( F_{0}+F_{\alpha}).
$$

\end{theorem}

After we prove the solvability of model equations estimate
\eqref{4.15.3} will allow us to prove the following existence
theorem.

\begin{theorem}
                                              \label{theorem 2.7.1}
There exists a unique $u\in \cC^{2+\alpha}(T)$ which satisfies
\eqref{2.6.4}.

\end{theorem}

   Theorem \ref{theorem 2.7.1} allows us to treat elliptic
equations.
\begin{theorem}
                                       \label{5.7.2}
Assume that the coefficients of $L$ and $f$ are independent of
$t$. Then there exists a unique $u\in C^{2+\alpha} $ satisfying $
Lu=f$.

\end{theorem}

   Proof. Let $u\in\cC^{2+\alpha}(0)$ be a unique solution
of the equation
\begin{equation}
                                                   \label{5.7.1}
u_{t}+Lu =f
\end{equation}
in $\bR^{d+1}_{0}$.  Since the coefficients and $f$ are
independent of $t$, for any fixed $s\geq0$, the function
   $u(s+\cdot , \cdot)$
   also satisfies \eqref{5.7.1} in $\bR^{d+1}_{0}$.
Obviously $u(s+\cdot,\cdot) \in\cC^{2+\alpha}(0)$. By uniqueness
$u(s+t,x)=u(t,x)$ for any $t\geq0$ and $x\in\bR^{d}$. In
particular, $u(s,x)=u(0,x)$,   $s \ge 0$,  and equation
\eqref{5.7.1} becomes $Lu=f$. This gives the existence of
solution. Uniqueness follows from assertion (ii) in   Theorem
\ref{theorem 1.10.1}. The theorem is proved.

   \begin{remark} \label{new}
   As we have pointed out in the Introduction,    a remarkable
feature of this theorem is that the set $LC^{2+\alpha}$ is
generally obviously wider than $C^{\alpha}$ unlike in the case
that the coefficients of $L$ are bounded   when
$LC^{2+\alpha}=C^{\alpha}$ always. Also notice that generally if
$c$ is bounded, the sets $L^{-1}C^{\alpha}$ are {\em different\/}
for different $L$ with growing coefficients because for a solution
$u$ of $Lu=f$ we have that $b^{i}D_{i}u$ is bounded. In addition,
the boundedness of $b^{i}D_{i}u$ means that the projection of the
gradient
    $Du(x)$  on $b(x)$    becomes smaller and smaller
   at points where $|b(x)|$ becomes large.
    What causes this
remains a mystery. It is certainly not true that this happens only
because $|Du(x)|$ is small where $|b(x)|$ is large, which is seen
if we take $a^{ij}=\delta^{ij}$, $c=1$, radially symmetric $f$,
and {\em any\/} $b^{i}(x)$ such that $b^{i}(x)x^{i}\equiv0$.
\end{remark}

   Our next main result concerns the Cauchy problem. Take and fix a
finite $S>T$ and introduce the space $\cC^{2+\alpha}(T,S)$ as the
set of functions $u\in\cC^{2+\alpha}(T)$ such that $u(t,x)=u(S,x)$
for $t\geq S$.

   We will consider the equation
\begin{equation}
                                           \label{4.15.12}
    u_{t}+Lu=f  \;\;\; \mbox {in $(T,S)\times\bR^{d}$.}
\end{equation}
Clearly, we may assume that $f= 0$  in
$(\bR\setminus[T,S])\times\bR^{d}$.   Therefore, concerning $f$ it
is enough to assume that
$$
|f(t,x)|\leq F_{0}\, c(t,x),\quad |f(t,x)-f(t,y)|\leq
F_{\alpha}|x-y|^{\alpha}, $$  whenever
    finite  $t\in (T,S)$
   and $x,y\in\bR^{d}$ are such that $|x-y|\leq1$.

\begin{theorem}
                                       \label{theorem 4.15.3}
Assume that we are given a function $g\in C^{2+\alpha} $.
Then there exists a unique $u\in\cC^{2+\alpha}(T,S)$ satisfying
   \eqref{4.15.12}   and such that $u(S,x)=g(x)$.
   Moreover, for this solution,
    for any finite  $t\in [T,S]$,
we have
\begin{equation}
                                           \label{4.15.10}
\|u(t,\cdot)\|_{2+\alpha} \leq
N(F_{0}+F_{\alpha}+\|g\|_{2+\alpha}),
\end{equation}
where $N=N(\delta,\alpha,d,K)$.
\end{theorem}

Proof. Uniqueness obviously follows from Theorem \ref{theorem
1.10.1}. To prove existence, introduce $\tilde{f}(t,x)=f(t,x)$ for
$t\leq S$ and
$$
\tilde{f}(t,x)=(\Delta+\partial/\partial t-\delta)
[e^{S-t}g(x)]=e^{S-t}[\Delta g(x)-(1+\delta)g(x)],
$$
for $t>S$. Also introduce the operator $\tilde{L}$ such that it
coincides with $L$ for $t\leq S$ and with
$\Delta+\partial/\partial t-\delta$ for $t>S$.

Then by Theorem \ref{theorem 2.7.1} there is a unique
$\tilde{u}\in\cC^{2+\alpha}(T)$ such that
\begin{equation}
                                           \label{4.15.11}
\tilde{u}_{t}+\tilde{L}\tilde{u}=\tilde{f}
\end{equation}
in $\bR^{d+1}_{T}$. Obviously $e^{S-t}g(x)$ is of class
$\cC^{2+\alpha}(S)$ and satisfies
\eqref{4.15.11} in $\bR^{d+1}_{S}$. By Theorem
\ref{theorem 2.7.1} we conclude that $\tilde{u}(t,x)=e^{S-t}g(x)$
for $t\geq  S$. In particular, $\tilde{u}(S,x)=g(x)$.

Also, $\tilde{u}$ satisfies \eqref{4.15.12}.  Now we define
$u(t,x)=\tilde{u}(t,x)$ for $t\leq S$ and
$u(t,x)=\tilde{u}(S,x)=g(x)$ for $t>S$. Then obviously
$u\in\cC^{2+\alpha}(T,S)$ and $u$ satisfies \eqref{4.15.12}.
Finally, estimate \eqref{4.15.10} follows immediately from
\eqref{4.15.3}, applied to $\tilde{u}$ and $\tilde{f}$, and the
definition of $\tilde{f}$. The theorem is proved.

\begin{remark}
                                        \label{remark 6.13.1}

   The previous theorem  can be adapted to the case
that    $c (t,x)$
     is only  a {\it  nonnegative} function
     on $\bR^{d+1}_T$
   with the other
   assumptions in Theorem \ref{theorem 4.15.3}
unchanged apart from \eqref{5.7.3} in Hypothesis \ref{H2}
   where we
replace $c$ with 1.

     Indeed,  if we want to solve
    equation  \eqref{4.15.12}
   with $c (t,x) \ge 0$ and
    final condition $g$,
       we can first
   solve the Cauchy problem
    \begin{equation} \label{di}
u_{t}(t,x)+Lu(t,x) -  u(t,x) =e^{ t-S }f(t,x)
   \end{equation}
in $(T,S)\times\bR^{d}$ with $u(S,x)=g(x)$,
   and then define $v (t,x) = e^{ S- t }
   u(t,x)$,
for $t \le S$ and $v(t,x) = g(x),$ $t>S $. Clearly,  $v
\in\cC^{2+\alpha}(T',S)$ for any finite $T'\in[T,S)$, $v(S, x) =
g(x)$, and $v$ solves \eqref{4.15.12}.
   By Theorem \ref{theorem 4.15.3}, we get,
    for any  finite $t\in[T,S]$,
$$
\| u(t, \cdot)\|_{2 + \alpha} \le N(
F_{0}+F_{\alpha}+\|g\|_{2+\alpha}),
$$
$$
\| v(t, \cdot)\|_{2 + \alpha} \le e^{S-t}
N(F_{0}+F_{\alpha}+\|g\|_{2+\alpha}).
$$
    Uniqueness for  
     \eqref{4.15.12} with $c (t,x) \ge 0$
      follows easily  from  the uniqueness already proved
       for equation \eqref{di}.

\end{remark}

\begin{remark} \label{cauchy}
   When  the operator $L$ has  coefficients independent of $t$ and
   $f =0$,
   Theorem \ref{theorem
   4.15.3}, in an obvious way (considering $S=0$, $T = - \infty$
    and inverting  time)
   allows one to introduce the
corresponding diffusion semigroup $T_{t}$ of bounded operators
mapping $C^{2+\alpha} $ into $C^{2+\alpha} $. By the
maximum principle we have
$$
\|T_{t}g\|_{0}\leq\|g\|_{0},\quad t \ge 0 ,\quad  g \in C^{2+
\alpha} ,
$$
and so an
   approximation argument allows one to extend $T_{t}$ from mappings
$C^{2+\alpha}  \to C^{2+\alpha} $ to mappings
$UCB(\bR^d)\to UCB(\bR^d)$ (where $UCB(\bR^d)$ stands for the
Banach space of all real uniformly continuous and bounded
functions defined  on $\bR^d$, endowed with the supremum norm).
Moreover, by interpolation theorems, $T_t$ will form a semigroup
of bounded operators mapping $C^{\alpha} $ into
$C^{\alpha} $.

   Several properties of the  diffusion
   semigroup $T_t$, corresponding to an operator
   $L$ with possibly unbounded
   time-independent coefficients
   $a$, $b$ and $c$ are investigated both  from an analytic point of
view (see  \cite{BL}
   \cite{CV}, \cite{DL}, \cite{Lu1}, \cite{LV},  \cite{MPW}) and from a
   probabilistic point of view  (see \cite{Ce0}, \cite{Ce}, \cite{DG},
   and the references therein).
\end{remark}

\mysection{Schauder estimates for equations with coefficients
independent of $x$}
                                           \label{section 4.16.2}

   In this section we concentrate on equations with the operator
$$
L_{0} u=a^{ij}(t)u_{x^{i}x^{j}},
$$
where $(a^{ij}(t))$ is a symmetric   matrix depending only on $t
\in (T, \infty)$ in a measurable way and such that
$$
K (\delta^{ij}) \geq(a^{ij} (t))\geq\delta (\delta^{ij}),
 \quad t > T,
$$

   For $t>s$ set
$$
A_{st}:=\int_{s}^{t}a(r)\,dr,\quad B_{st}:=A_{st}^{-1}.
$$
Observe that the matrices $A_{st}$ are  nondegenerate so that
$B_{st}$ is well defined and
$$
K(t-s)|\xi|^{2}\geq A_{st}^{ij} \xi^{i}\xi^{j}
\geq\delta(t-s)|\xi|^{2},
$$
$$
\delta^{-1}(t-s)^{-1}|\xi|^{2}\geq B_{st}^{ij} \xi^{i}\xi^{j} \geq
K^{-1}(t-s)^{-1}|\xi|^{2},\;\;\; t>s,\;\; \xi \in \bR^d.
$$
Define ($I_{t>s}$ stands for the indicator function of the set
$\{t>s\}$)
\begin{equation} \label{pp}
p(s, t,x)=I_{t>s}(4\pi)^{-d/2}
(\text{det}\,B_{st})^{1/2}\exp(-(B_{st}x,x)/4),
\end{equation}
$$
G_{}f(s,x)=\int_{s}^{\infty}\int_{\bR^{d}} p(s,t,x-y)f(t,y) \,dydt
$$
$$
= \int_{0}^{\infty}\int_{\bR^{d}} p(s,r+s,x-y)f(r+s,y) \,dydr,
$$
for any bounded measurable function $f$  on $\bR^{d+1}$  with
compact support.

   The next  two lemmas
    concerning the potential $G$ are
essentially known even if in the literature they are stated in
    a slightly different form (see \cite{B}, \cite{Kn}, \cite{Kry2},
\cite{L}, \cite{Lo}). For the sake of completeness we
   still give   the  proofs albeit sometimes somewhat
   sketchy.

\begin{lemma}
                                         \label{lemma 4.24.1}
If $u\in\cC^{2+\alpha}_{0}(T)$, then for any  finite $t\geq T$ and
$x\in\bR^{d}$ we have
\begin{equation}
                                                  \label{4.24.1}
u(t,x)=-G (  u_{t}+L_{0}u)(t,x).
\end{equation}

\end{lemma}

Proof. Set $f = u_t + L_0 u $ and observe that
   $f$ is a bounded
measurable function defined in $\bR^{d+1}_{T}$ 
vanishing for large $|t|+|x|$. We
have, for any finite $t>s\geq T$, $x \in \bR^d$,
$$
u(t,x) - u(s,x) = \int_s^t f(r,x)dr - \int_s^t L_0 u(r, x) dr.
$$

 Using the Fourier transform in the space variable $x$ we get
$$
\hat u(t,\xi) - \hat u(s,\xi) = \int_s^t \hat f(r,\xi)dr +
\int_s^t \xi^i \xi^j a^{ij}(r) \hat u(r, \xi) dr,
$$
i.e., fixing $\xi \in \bR^d$, we have a.s. in $t$, $\hat u_t(t,
\xi) = \hat f (t, \xi) + \xi^i \xi^j a^{ij}(t) \hat u(t, \xi)$. It
follows that
$$
\hat u (t, \xi) = - \int_t^{ \infty} e^{- {A^{ij}_{tr}\xi^i \xi^j
}} \hat f (r, \xi) dr.
$$
Taking the anti-Fourier transform we easily get the assertion.  
\begin{lemma}
                                         \label{lemma 4.24.2}
If $f(t,x)$ is a bounded  measurable function in $\bR^{d+1}$
vanishing if $t\geq S$, for some 
(finite) constant $S$, and such that
\begin{equation}
                                                  \label{4.13.2}
\sup_{t>T}[f(t,\cdot)]_{\alpha}<\infty,
\end{equation}
then, for each finite $t\geq T$, the function $Gf(t,x)$ is twice
continuously differentiable in $x$ and, for any $x,y\in\bR^{d},$
   we have
\begin{equation}
                                                  \label{4.13.1}
|D^{2}G f(t,x) |\leq N  \sup_{s>t}[f(s,\cdot)]_{\alpha},
\end{equation}
where $N$ depends only on $(S-t)_{+},d,\delta$, and $\alpha$, and
\begin{equation}
                                                  \label{4.24.01}
|D^{2}G f(t,x)-D^{2}G f(t,y)|\leq N
   \sup_{s>t} [f(s,\cdot)]_{\alpha}|x-y|^{\alpha},
\end{equation}
where $N$ depends only on $d,\delta$, and $\alpha$.
\end{lemma}

    Proof.
    Apart from the fact that the constants $N$
    are
   independent of $K$, this   is   a
   standard   result  by now  (it was first proved in  \cite{B}
    when $a^{ij}$ are independent on $t$).
   It can  also be
   proved by adapting  the computations  in Section IV.2 of \cite{LSU}
   (see, in particular, pages 276 and 277 where the case of
    $a^{ij}$  independent on $t$ is treated) or, in a more direct
    way,   arguing as in the proof
      of  Theorem 4.2 in \cite{P1}.

   Estimates \eqref{4.13.1} and \eqref{4.24.01} are proved
   in Lemma 3.2 of \cite{Kry2}  for measurable
$a^{ij}$  in a slightly sharper form (with the maximal function of
$[f(s,\cdot)]_{\alpha}$ in place of $\sup_{s> t}$, which allowed
one to consider spaces with norms that are $C^{\alpha}$ in $x$ and
$L_{p}$ in $t$). However, $f$ in \cite{Kry2} was assumed to be in
$C^{\infty}_{0} (\bR^{d+1})$.

Actually, in the proof of Lemma 3.2 of \cite{Kry2} the facts that
$f$ is so  much regular in $t$ and vanishes for large $|x|$ were
never used and what was used is that $f$ has compact support in
$t$, is bounded, and the derivative of $f(s,x)$ in $x$ up to the
third order are continuous in $x$ and bounded in $(s,x)$. To relax
further these requirements and include the class of $f$ under
consideration we introduce
$f^{(\varepsilon)}(t,x)=(f(t,\cdot)*\zeta_{\varepsilon})(x)$,
where
$\zeta_{\varepsilon}(x)=\varepsilon^{-d}\zeta(x/\varepsilon)$ and
$\zeta$ is a nonnegative $C^{\infty}_{0}$ function on $\bR^{d}$
which integrates to one. Owing to \eqref{4.13.2} we have
$$
\|f-f^{(\varepsilon)}\|_{0,T}\leq N\varepsilon^{\alpha},\quad
\|Gf-Gf^{(\varepsilon)}\|_{0,T}\leq N\varepsilon^{\alpha},
$$
$$
\sup_{s>t}[f^{(\varepsilon)}(s,\cdot)]_{\alpha} \leq \sup_{s>t}
[f(s,\cdot)]_{\alpha}.
$$
Hence, for each $t$, as $\varepsilon\to0$, the functions
$Gf^{(\varepsilon)}(t,\cdot)$ converge to $Gf(t,\cdot)$ uniformly
and, owing to \eqref{4.13.1} and \eqref{4.24.01}, have uniformly
bounded and uniformly continuous second-order derivatives in $x$.
A standard result from Calculus implies that, for each $t$,
$Gf(t,x)$ is twice continuously differentiable in $x$ and
$D^{2}Gf^{(\varepsilon)}(t,\cdot) \to D^{2}Gf (t,\cdot)$. By
passing to the limit in \eqref{4.13.1} and \eqref{4.24.01} with
$f^{(\varepsilon)}$ in place of $f$ we obtain the desired result.
This will again prove our estimates with constants depending on
$K$. The fact, which will not be used in this article,
   that the constants $N$ are independent on $K$
follows from a general result proved in \cite{Kry3}.

   We  need
   a result on the unique solvability of the heat
   equation  for which we could not find a precise
    reference
    in  the literature.  A slight difficulty
in proving it is caused by  the fact that
   the datum $f$ is only measurable in time (indeed, if $f$ is
continuous
    also in time, the  proof is straightforward  and the result
    becomes well known).

   \begin{lemma}
                                        \label{lemma 4.14.1}
Let $T>-\infty$ and assume that $f(t,x)$ is a bounded measurable
function vanishing for $t\geq S$, where $S$ is a constant, and
such that \eqref{4.13.2} holds. Then there exists a unique
function $u\in\cC^{2+\alpha}(T)$ such that in $\bR^{d+1}_{T}$ we
have
\begin{equation}
                                                 \label{4.14.5}
u_{t}+\Delta u-\delta u=f.
\end{equation}
Furthermore, $u(t,x)=0$ for $t\geq S$.
\end{lemma}

Proof. Uniqueness follows from Theorem \ref{theorem 4.25.1}, which
we prove later in a much more general setting.

   To prove existence set $g(t,x)=-e^{\delta t}f(t,x)$
and make the change of the unknown function $v(t,x)=e^{\delta
t}u(t,x)$. Then \eqref{4.14.5} becomes
\begin{equation}
                                                 \label{4.14.6}
v_{t}+\Delta v =-g.
\end{equation}

Lemma \ref{lemma 4.24.1} suggests a natural candidate for $v$:
$$
v(t,x)=G_{0}g(t,x),
$$
where by $G_{0}$ we denote the operator $G$ constructed from
$a^{ij}\equiv\delta^{ij}$. Therefore, we introduce $v$ by the
above formula and proceed further.

Obviously $v(t,x)$
   is bounded and vanishes for $t\geq S$.
Owing to Lemma \ref{lemma 4.24.2} and interpolation inequalities,
$v\in\cC^{2+\alpha}(T)$. By the way, this and the fact that
$T>-\infty$ imply that $u(t,x):=e^{-\delta t}v(t,x)$ is of class
$\cC^{2+\alpha}(T)$ as well. Now the relation between
\eqref{4.14.5} and \eqref{4.14.6} shows that to finish proving the
lemma it only remains to prove that \eqref{4.14.6} holds.

First, observe that in our case $p(t,s,x)=p(0,s-t,x)=:p(t-s,x)$
   (see \eqref{pp})
and define the heat semigroup by
$$
T_{t}h(x)=(p(t,\cdot)*h)(x),\quad t>0
$$
acting   on bounded measurable   functions $h$ defined on
$\bR^{d}$. Then, take an $\varepsilon>0$ and introduce
$g^{(\varepsilon)}(t,\cdot) =T_{\varepsilon}g(t,\cdot)$,
\begin{equation}
                                                      \label{4.15.1}
v^{(\varepsilon)}(t,x)=G_{0}g^{(\varepsilon)}   (t,x)
=\int_{t}^{S}T_{s-t}T_{\varepsilon}g(s,x)\,ds
=\int_{t}^{S}T_{s-t+\varepsilon} g(s,x)\,ds.
\end{equation}
One knows that, for each $x$ and bounded $h$, the function
$T_{\tau} h( x)$ as a function of $\tau$ is continuously
differentiable in $\tau$ for $\tau>0$, infinitely differentiable
in $x$, and
$$
\frac{\partial}{\partial \tau}T_{\tau}h(x) =\Delta T_{\tau} h(x).
$$
In particular, for any $\tau$ the rules of differentiating
integrals depending on parameters imply that
$$
\frac{\partial}{\partial s}T_{s-t+\varepsilon} g(\tau,x) =\Delta
T_{s-t}T_{ \varepsilon} g(\tau,x) = T_{s-t}\Delta
g^{(\varepsilon)}(\tau,x)
$$
for $s>t$. It follows that
$$
T_{s-t+\varepsilon} g(\tau,x)=T_{\varepsilon}g(\tau,x)
+\int_{t}^{s}\frac{\partial}{\partial r}T_{r-t+\varepsilon}
g(\tau,x) \,dr
$$
$$
=T_{\varepsilon}g(\tau,x) +\int_{t}^{s}T_{r-t}\Delta
g^{(\varepsilon)}(\tau,x)\,dr.
$$
We plug in here $s$ in place of $\tau$ and go back to
\eqref{4.15.1}. Then we find
\begin{equation}
                                                      \label{4.15.2}
v^{(\varepsilon)}(t,x)=\int_{t}^{S} g^{(\varepsilon)}(s,x)\,ds
+\int_{t}^{S}\int_{t}^{s}T_{r-t}\Delta
g^{(\varepsilon)}(s,x)\,drds.
\end{equation}
Changing the variable $r$ by $\tau=s-r+t$ and then using Fubini's
theorem   we obtain
$$
\int_{t}^{S}\int_{t}^{s}T_{r-t}\Delta g^{(\varepsilon)}(s,x)\,drds
=\int_{t}^{S}\int_{t}^{s}T_{s-\tau}\Delta
g^{(\varepsilon)}(s,x)\,d\tau ds
$$
$$
= \int_{t}^{S} \big(\int_{\tau}^{S}T_{s-\tau} \Delta
g^{(\varepsilon)}(s,x)\,ds \big)d\tau= \int_{t}^{S}\Delta
v^{(\varepsilon)}(\tau,x)\,d\tau.
$$
Upon combining this with \eqref{4.15.2} we conclude that
$v^{(\varepsilon)}$ satisfies \eqref{4.14.6} with
$g^{(\varepsilon)}$ in place of $g$.

Finally, we send $\varepsilon\downarrow0$ and observe that
$g^{(\varepsilon)}\to g$ uniformly in $\bR^{d+1}_{T}$ since $g$ is
uniformly continuous in $x$. It follows that $v^{(\varepsilon)}\to
v$ uniformly in $\bR^{d+1}_{T}$. Furthermore,
   $D^{2}G_{0}g^{(\varepsilon)}\to D^{2}G_{0}g$
uniformly in $\bR^{d+1}_{T}$ since for any $\beta\in(0,\alpha)$,
by \eqref{4.13.1}, we have
$$
\|D^{2}G_{0}(g^{(\varepsilon)}-g)\|_{0, T}  \leq N \sup_{s>
T}[g^{(\varepsilon)}(s,\cdot)-g(s,\cdot)]_{\beta} \to0,
$$
   as $\varepsilon\downarrow0$,  where the last
relation follows from the fact that \eqref{4.13.2} holds with $g$
in place of $f$  (just in case,
   note that $N$ depends  on
   $(S-T)_+$).
   This argument allows us to pass to the limit in the
integral version of \eqref{4.14.6}. The lemma is proved.

 The next elementary result is well known.
\begin{lemma}
                                      \label{lemma 6.13.1}
Let $\gamma\geq1$ and let
  $Q$  be a convex closed round   cone in $\bR^d$
   with vertex at the origin such that
for any unit ball which lie inside $Q$ the distance of its center
to the origin is greater than or equal to $\gamma$. Let $(u^{ij})$
be a $d\times d$ symmetric matrix. Then there is a constant
$N=N(\gamma,d)$ such that for any $i,j=1,...,d$ we have
$$
|u^{ij}|\leq N\max_{|\xi|=1,\xi\in Q}\big|
\sum_{i,j=1}^{d}u^{ij}\xi^{i}\xi^{j}\big|.
$$
\end{lemma}

   Now follows an embedding result
generalizing
   \cite[Lemma 1]{Kn}. Its proof uses parabolic dilations
    and  is simpler   than  the  one in \cite{Kn}.
     Its generality is
   in part
motivated by possible applications to boundary value
    problems.

Take $\gamma$ and  $Q$   as in Lemma \ref{lemma
6.13.1}, take an $h>0$, and consider the truncated cone
   $$
Q_{h}= Q\cap \{x:|x|\leq h\} .
    $$
    The  spaces   $C^{\alpha} (Q_{h})$
    are defined in the same way as $C^{\alpha}=C^{\alpha}
(\bR^d)$. We also write $[\cdot]_{\alpha, Q_{h}}$ to denote
     the usual H\"older seminorm in $C^{\alpha} (Q_{h})$.
Similarly we introduce the spaces $C^{2+ \alpha}(Q_{h})$ and the
seminorms $[\cdot]_{2+\alpha, Q_{h}}$. Notice that $Q_{h}$ is a
closed set. In particular, the functions from $C^{2+
\alpha}(Q_{h})$ are twice continuously differentiable in the
interior of $Q_{h}$ and their derivatives admit continuous
extension to the boundary of $Q_{h}$. In the following lemma by
$D^{2}u(r,0)$ and $Du(r,0)$ we mean these continuations.

\begin{lemma}

                                        \label{lemma 4.13.2}
     Let $u: [0,h^{2}] \times Q_{h} \to \bR$ be a
   continuous  function such that
    $u(t, \cdot ) \in C^{2+ \alpha}  (Q_{h})$, $t \in [0,h^{2}]$,
   and assume that there
   exists a   function
    $g(t,x)$ defined on $[0,h^{2}] \times Q_{h}$
    such that  $ g(t, \cdot ) \in C^{\alpha } (Q_{h})$,
$t \in [0,h^{2}]$, and
   \eqref{g1} is verified for $0 \le s\le t \le h^{2} $,
$x \in Q_{h}$
   (we set $u_t =g $).

Then there is a constant $N=N(\gamma,d)$ such that
\begin{equation}
                                                  \label{4.13.3}
|D^{2}u(h^{2},0)-D^{2}u(0,0)|\leq N I_{h}h^{\alpha},
  \end{equation}
\begin{equation}
                                                  \label{4.14.1}
|Du(h^{2},0)-Du(0,0)|\leq NI_{h}h^{1+\alpha},
  \end{equation}
where
$$
  I_{h}:=\sup_{r\in[0,h^{2}]} \big([u_{t}(r,\cdot)]_{\alpha, Q_{h}}
+[D^2 u (r,\cdot)]_{ \alpha ,Q_{h} }\big).
$$

\end{lemma}
   Proof.  First we deal with \eqref{4.13.3}.
The parabolic dilation  $(s,x) \mapsto (4^{-1}h^{-2}s,2^{-1}h^{-1}
x )$  allows us to assume that $h=2$.

  Next,   assume that the first basis vector $\ell$
is inside $Q_2$ ($=Q_{h}$). Write
$$
|D_{11} u(4,0)-D_{11}u(0,0)|   \le   |D_{11} u(4,0)-[u(4,
2\ell)-2u(4,\ell)+u(4,0)   ] |
$$
$$
+|D_{11} u(0,0)-[u(0 , 2\ell)-2u(0,\ell)+u(0,0) ] |+I_{+}+I_{-},
$$
where
$$
I_{\pm}=|[u(4, \ell \pm\ell)-u(0, \ell \pm\ell)]-
[u(4,\ell)-u(0,\ell)]|.
$$
Taylor's formula shows that
$$
|D_{11} u(4,0)-[u(4, 2\ell)-2u(4,\ell)+u(4,0)   ]  | \leq 2
     [D_{11}u(4,\cdot)]_{ \alpha , Q_{2} },
$$
$$
   |D_{11} u(0,0)-[u(0 , 2\ell)-2u(0,\ell)+u(0,0) ]|
  \leq 2
     [D_{11}u(0,\cdot)]_{ \alpha , Q_{2} }.
$$
By the Newton-Leibnitz formula
$$
I_{\pm} =\int_{0}^{4}[u_{t}(r,\ell \pm\ell)-u_{t}(r, \ell )]\,dt,
$$ which implies
that
$$
I_{\pm}\leq\sup_{r\in[0,4]}[u_{t}(r,\cdot)]_{{{\alpha}, Q_{2}}}.
$$
Upon combining the above estimates,  we come to
$$
|D_{11}u(4,0)-D_{11}u(0,0)|\leq 4I_{2}.
$$

Having this estimate proved for  $\xi^{i}\xi^{j}D_{ij}u$ with
$\xi$ being the first basis vector, under the assumption that it
lies inside $Q_{2}$, we also have
  $$
|\xi^{i}\xi^{j}D_{ij}u(4,0)-\xi^{i}\xi^{j}D_{ij}u(0,0)| \leq
4I_{2}
$$
for all unit $\xi\in Q_{2}$. Applying Lemma \ref{lemma 6.13.1}
yields \eqref{4.13.3}.

In light of \eqref{4.13.3},
estimate \eqref{4.14.1} is, actually, classical and we give its
proof just for completeness. 
Again we may assume that $h=2$ and first consider the case when
$\ell\in Q_{2}$. Then
$$
|D_{1}u(4,0)-D_{1}u(0,0)|\leq |J_{1}(1)|+|J_{2}|,
$$
where
$$
J_{1}(r)=[u(4,r\ell)-u(4,0)-rD_{1}u(4,0)]
-[u(0,r\ell)-u(0,0)-rD_{1}u(0,0)],
$$
$$
J_{2}=[u(4, \ell)-u(4,0)]-[u(0, \ell)-u(0,0)]=\int_{0}^{4}
[u_{t}(r,\ell)-u_{t}(r,0)]\,dr.
$$
We estimate $|J_{2}|$ in the same way as $I_{\pm}$ above and
notice that by Taylor's formula (or by integrating by parts)
$$
J_{1}(r)=\int_{0}^{r}(r-z)[D_{11}u(4,z\ell)-
D_{11}u(0,z\ell)]\,dz.
$$
We can certainly take the point $z\ell$ as the origin in
$\bR^{d}$ and use \eqref{4.13.3}. Then we see that $|J_{1}(1)|\leq
NI$ and this leads to the estimate
$$
|\xi^{i}D_{i}u(4,0)-\xi^{i}D_{i}u(0,0)|\leq NI
$$
if $\xi=\ell$. The same estimate holds for any unit vector $\xi\in
Q_{2}$ and this implies \eqref{4.14.1}. The lemma is proved.

By shifting the origin in $\bR^{d+1}$ to points $(s,x)$ and
denoting $t-s=h^{2}$ we obtain the following.

\begin{corollary}

                                        \label{corollary 4.13.2}
For any  $u\in\cC^{2+\alpha}(T)$ and finite $t,s\geq T$ and
$x\in\bR^{d}$ we have
\begin{equation}
                                             \label{6.14.3}
|D^{2}u(t,x)-D^{2}u(s,x)|\leq N I|t-s|^{\alpha/2},
\end{equation}
$$
|D u(t,x)-D u(s,x)|\leq N I|t-s|^{(1+\alpha)/2},
$$
where $N$ depends only on $d$ and
$$
I=\sup_{r\in[s,t]} \big([u_{t}(r,\cdot)]_{\alpha} +[u
(r,\cdot)]_{2+\alpha}\big).
$$
\end{corollary}

\begin{remark}
                                              \label{remark do}
   In connection with  Corollary \ref{corollary 4.13.2}
   it is also worth
noting that  if $u\in\cC^{2+\alpha}(T)$, then $u$ is locally
Lipschitz in $(t,x)\in\bR^{d+1}_{T}$ since the derivative
   $D u$  is bounded and $u_{t}$ is locally bounded.
\end{remark}

\begin{lemma}
                                        \label{lemma 2.6.1}
If $u\in\cC^{2+\alpha}_{0}(T)$, then for any $x\in\bR^{d}$ and
finite $t,s\ge  T$   we have
\begin{equation}
                                                  \label{1.10.1}
[u(t,\cdot)]_{ 2+\alpha}\leq N \sup_{r > t}[(u_{t}+
L_{0}u)(r,\cdot)]_{\alpha},
\end{equation}
\begin{equation}
                                                  \label{2.6.1}
|D^{2}u(s,x)-D^{2}u(t,x)|\leq N  \sup_{r > t}[(u_{t}+
L_{0}u)(r,\cdot)]_{\alpha}|s-t|^{\alpha/2},
\end{equation}
where the constant $N$ depends only on $d,\delta,\alpha$ (and is
independent of~$K$).
\end{lemma}
    To prove the lemma it suffices to make
just few comments. Indeed, \eqref{1.10.1} follows directly from
Lemmas \ref{lemma 4.24.1} and \ref{lemma 4.24.2}. Estimate
\eqref{2.6.1} with $f$ in place of $u_{t}+L_{0}u$ follows from
Theorem 3.3 of \cite{Kry2} for $u$ in the form $Gf$ if $f\in
C^{\infty}_{0} (\bR^{d+1})  $. By Lemma \ref{lemma 4.24.1} any   $
u \in \cC^{2+\alpha}_{0}(T)$  has this form with $f$ having less
regularity than it is required in \cite{Kry2}.
   Then one can repeat what is said concerning the proof of
   Lemma \ref{lemma 4.24.2}.

\begin{remark}
It is also worth noting that \eqref{2.6.1} will not be used in the
future and if one does not care about the statement that $N$ in
\eqref{2.6.1} is independent of $K$, then one can get the estimate
from   Corollary \ref{corollary 4.13.2} and \eqref {1.10.1} since
$$
[u_{t}(r,\cdot)]_{\alpha}\leq [(u_{t}+ L_{0}u)(r,\cdot)]_{\alpha}+
[ L_{0}u (r,\cdot)]_{\alpha}
$$
$$
\leq [(u_{t}+ L_{0}u)(r,\cdot)]_{\alpha}+ N
[u(r,\cdot)]_{2+\alpha},
$$
   where $N =N (K, d)$.
   One can also note that, in turn, \eqref{6.14.3}
follows from \eqref{2.6.1}, which is seen if one takes
$L_{0}=\Delta$.
\end{remark}

   Here is a rather surprising generalization of estimate
\eqref{1.10.1}.

\begin{lemma}
                                            \label{lemma 1.10.1}
Let $N$ be the constant from \eqref{1.10.1}.

(i) If $u\in\cC^{2+\alpha}(T)$ and there is an $R\in(0,\infty)$,
such that $u(t,x)=0$ whenever $|x|\geq R$ and finite $t\ge  T$, then
for any finite $t\ge  T$, any locally bounded, nonnegative,  measurable
function $c_{0}(t)$, and any locally bounded, $\bR^{d}$-valued,
measurable function $b_{0}(t)$ we have
\begin{equation}
                                                 \label{first}
[u(t,\cdot)]_{ 2+\alpha} \leq N \sup_{s> t}[(u_{t}+
L_{0}u+b^{i}_{0}u_{x^{i}}-c_{0}u) (s,\cdot)]_{\alpha}.
\end{equation}

   (ii) If the coefficients of $L_{0}$ are independent of $t$,
then for any $u\in C^{2+\alpha} $ with compact support,
any constant $c_{0}\geq0$, and any constant vector
$b_0\in\bR^{d}$, we have
$$
[u]_{ 2+\alpha}\leq N [L_{0}u+b^{i}_{0}u_{x^{i}}-c_{0}u]_{\alpha}.
$$
\end{lemma}
Proof.  Assertion (ii) follows directly from (i) if we take in the
latter $u$ independent of $t$.

We prove (i) first assuming that $u\in\cC^{2+\alpha}_{0}(T)$. For
$t\in\bR$ define
$$
B(t)=\int_{0}^{t}b_{0}(s)\,ds,\quad v(t,x)=u(t,x+B(t)).
$$
Since $u\in\cC^{2+\alpha}_{0}(T)$ and the derivative of $B$ is
locally bounded, $v\in\cC^{2+\alpha}_{0}(T)$. By plugging in $v$
in place of $u$ in \eqref{1.10.1} we obtain \eqref{first} if
$c_{0}\equiv0$.

To let $c_{0}$ enter into the picture introduce
$$
C(t)=\int_{0}^{t}c_{0}(s)\,ds,\quad v(t,x)=e^{-C(t)}u(t,x).
$$
Again $v\in\cC^{2+\alpha}_{0}(T)$ and by substituting $v$ in place
of $u$ in \eqref{first} with $c_{0}\equiv0$ we get
$$
[u(t,\cdot)]_{ 2+\alpha}e^{-C(t)} \leq N \sup_{s> t}[(u_{t}+
L_{0}u+b^{i}_{0}u_{x^{i}}-c_{0}u) (s,\cdot)]_{\alpha}e^{-C(s)},
$$
$$
[u(t,\cdot)]_{ 2+\alpha} \leq N \sup_{s> t}[(u_{t}+
L_{0}u+b^{i}_{0}u_{x^{i}}-c_{0}u)
(s,\cdot)]_{\alpha}e^{C(t)-C(s)}.
$$
In case $u\in\cC^{2+\alpha}_{0}(T)$, this yields \eqref{first} in
full generality since $C(s)\geq C(t)$ for $s>t$.

    To pass to the case of general  $u$, take
a $\zeta\in C^{\infty}_{0}(\bR)$, such that $\zeta(0)=1$ and
$0\leq\zeta\leq1$, and substitute
$$
u^{n}(s,x):=\zeta(s/n)u(s,x)
$$
in \eqref{first}. Then let $n\to\infty$. After that it will only
remain to observe that
$$
[(u^{n}_{t}+ L_{0}u^{n}+b^{i}_{0}u^{n}_{x^{i}}-c_{0}u^{n})
(s,\cdot)]_{\alpha}\leq n^{-1}  [u(s,\cdot)]_{\alpha}\sup|\zeta'|
$$
$$
+[u_t + L_{0}u+b^{i}_{0}u_{x^{i}}-c_{0}u) (s,\cdot) ]_{\alpha}.
$$
The lemma is proved. \mysection{Proof of the main results}

                                           \label{section 4.16.3}

   We start with proving the following a priori estimate
   of the kind of the maximum principle, in which the assumptions
on the coefficients are      weaker than Hypotheses \ref{hy} and
\ref{H2}.

   The proof of the next result seems to be new even in the case
of bounded coefficients $a,$ $b$, $c$ and bounded datum $f$.  It
is important for the future to observe that condition (iii) of
Theorem \ref{theorem 4.25.1} below is satisfied for any
   $u\in\cC^{2+\alpha}(T)$,
    which follows from  Corollary \ref{corollary 4.13.2}
and Remark~\ref{remark do} applied to $u(t,x)\zeta(x)$, where
$\zeta$ is any function of class $C^{\infty}_{0} (\bR^{d})$. It is
also worth noting that the result of Theorem \ref{theorem 4.25.1}
cannot be obtained from the Alexandrov maximum principle for
parabolic equations since it requires $u\in
W^{1,2}_{d+1,loc}(\bR^{d+1}_{T})$ and in the theorem $u_{t}$ may
be only locally summable to the first power. However, by using
parabolic Aleksandrov estimates one could considerably relax the
assumptions on $u$ if $c$ and $f$ are locally bounded.

\begin{theorem}
                                  \label{theorem 4.25.1}
(i) Assume that  in $\bR^{d+1}_{T}$ we have an operator $L$  as in
\eqref{lu} and $f$ such that   $a^{ij}$, $b^i$, $c$, $f$ are
measurable functions, $(a^{ij} (t,x)) $ is symmetric  and
nonnegative, $c(t,x) \geq\delta$, $-f (t,x)\leq F_{0}c(t,x)$,
$$
|a(t,x)|\leq K(t)(1+|x|^{2}), \quad |b(t,x) |\leq K(t)(1+|x|)
$$
   in $\bR^{d+1}_{T}$, where $K(t)$ is a locally bounded function
on $\bR$;

(ii) Assume that, for any $x$, the functions
 $c(t ,x)$ and $f(t,x)$ are
locally integrable   in $t$  on $\bR$;

(iii) Assume that in $\bR^{d+1}_{T}$ we are given a bounded
continuous function $u(t,x)$ which is twice continuously
differentiable in $x$ for each finite   $t\geq T$ and such that $ u_{x}$
and $u_{xx}$ are continuous with respect to $(t,x)$ in
$\bR^{d+1}_{T}$; 

(iv) Finally, assume that \eqref{2.7.5} holds for
each $x\in\bR^{d}$ and  finite
    $t>s$ such that $s\geq T$. Then in $\bR^{d+1}_{T}$
\begin{equation}
                                             \label{4.25.2}
u(t,x)\leq F_{0}.
\end{equation}

Furthermore, if $|f (t,x)| \leq F_{0}c(t,x)$, then $|u(t,x)|\leq
F_{0}$ in $\bR^{d+1}_{T}$.

\end{theorem}
Proof. Obviously, the second assertions follows from the first
one. To prove the first assertion observe that the function
$v=u-F_{0}$ satisfies \eqref{2.6.4} with
    $f$ replaced by
$g=cF_{0}+f\geq0$ and $-g\leq 0\cdot c$.  If our result is true
for $f\geq0$ and $F_{0}=0$, then $v\leq0$ and $u\leq F_{0}$. It
follows that in the rest of the proof we may confine ourselves to
the case that $F_{0}=0$. Without losing generality we may also
assume that $T>-\infty$ and even that $T=0$.

Next notice that to prove \eqref{4.25.2} (with $F_{0}=0$), it
suffices to prove that for any $S>0$, $t\in[0,S]$, and
$x\in\bR^{d}$ we have
\begin{equation}
                                             \label{4.26.2}
u(t,x)\leq e^{\delta'(t-S)}\sup_{y\in\bR^{d}}u_{+}(S,y),
\end{equation}
where $\delta'=\delta/2$. Indeed, one can then let $S\to\infty$
and use that  $u$ is bounded by assumption. In turn, to prove
\eqref{4.26.2} it suffices to show   that for any $\gamma>0$
\begin{equation}
                                             \label{4.26.1}
\bar{u}(t,x):=(u(t,x)-\gamma
   w(t,x)
   )e^{-\delta' t}\leq
e^{-\delta' S}\sup_{y\in\bR^{d}}u_{+}(S,y),
\end{equation}
where
$$
w(t,x)= (1+|x|^{2})e^{-N_{0}t}
$$
and we choose the constant $N_{0}$
   so large that in $[0,S]\times \bR^{d}$ we
have
$$
w_{t}+ Lw=:g  <0,
$$
which is clearly possible.

We use the fact that $u$ is bounded and continuous, $w$ tends to
infinity as
   $|x|\to\infty$ to conclude that   at certain point $(t_{0},x_{0})\in
[0,S]\times\bR^{d}$ the function $\bar{u}(t,x)$ takes its maximum
value over $[0,S]\times\bR^{d}$. If
$$
\bar{u}(t_{0},x_{0})\leq
    e^{-\delta' S}
\sup_{y\in\bR^{d}}u_{+}(S,y),
$$
then \eqref{4.26.1}  obviously holds. Therefore, it suffices to
show that the inequality
\begin{equation}
                                             \label{4.26.3}
\bar{u}(t_{0},x_{0})>     e^{-\delta' S}
\sup_{y\in\bR^{d}}u_{+}(S,y),
\end{equation}
is impossible.

We argue by contradiction and suppose   that \eqref{4.26.3} holds.
Then obviously $0\leq t_{0}<S$ and $\bar{u}(t_{0},x_{0})>0$.
Furthermore, at $(t_{0},x_{0})$ we have that $D\bar{u} =0$ and
the      Hessian   matrix $D^{2}\bar{u} $ is nonpositive, which
   implies that, for any $\tau>0$, there exists a
$\theta>0$ such that
$$
|D \bar{u} (t,x_{0})|\leq\tau,\quad D^{2}\bar{u} (t,x_{0})\leq\tau
(\delta^{ij}),\quad \bar{u} (t,x_{0})>0
$$
whenever $0\leq t-t_{0}\leq\theta$. Also
\begin{equation}
                                             \label{4.26.4}
0\geq \bar{u}(t,x)-\bar{u}(t_{0},x_{0})=
\int_{t_{0}}^{t}\bar{u}_{t}(s,x_{0})\,ds.
\end{equation}

However,
$$
\bar{u}_{t}=-\delta'\bar{u} +  (u_{t}-\gamma w_{t}) e^{-\delta't}
$$
$$
=-\delta' \bar{u} - L\bar{u}+fe^{-\delta't}-\gamma ge^{-\delta't}
\geq -\delta'\bar{u} - L\bar{u}
$$
and the last expression at points $(t,x_{0})$ such that $0\leq
t-t_{0} \leq\theta$ is greater than
$$
-N\tau+(c-\delta')\bar{u}\geq -N\tau+\delta'\bar{u},
$$
where the constant $N$ is independent of $t$ and $\tau$. Hence,
for   $ t_0 + \theta \ge t \geq t_{0}$, equation \eqref{4.26.4}
implies that
$$
0\geq-N\tau(t-t_{0})+\delta'\int_{t_{0}}^{t} \bar{u}(s,x_{0})\,ds.
$$
We divide both part of this inequality by $t-t_{0}$, let
$t\downarrow t_{0}$, and use the continuity of $\bar{u}$. Then we
conclude that $\delta'\bar{u}(t_{0},x_{0})\leq N\tau$ and since
$\tau>0$ is arbitrary, $\bar{u}(t_{0},x_{0})\leq0$, which yields
the desired contradiction with \eqref{4.26.3}. The theorem is
proved.

{\bf Proof of Theorem \ref{theorem 1.10.1}}. As usual the second
assertion is obtained from the first one by taking there $u$
independent of $t$.

   While proving (i), first observe that it suffices to prove the
estimate
\begin{equation}
                                                          \label{2.7.1}
\sup_{t\ge  T}[u(t,\cdot)]_{2+\alpha}\leq N  (F_{0}+F_{\alpha})
+N\sup_{t\ge  T}\|u(t,\cdot)\|_{2}.
\end{equation}
   Indeed, once \eqref{2.7.1} is proved,
\eqref{4.15.3} follows easily from the interpolation inequality
$$
\|v\|_{2}\leq N(\varepsilon)\|v\|_{0}+\varepsilon[v]_{2+\alpha},
\quad \varepsilon>0, \;\;\; v \in C^{2+ \alpha} ,
$$
   and the fact  that $|u(t,x)| \leq F_{0}$
   in $\bR^{d+1}_{T}$ by
Theorem \ref{theorem 4.25.1}.

To prove \eqref{2.7.1} fix an $\varepsilon\in(0,1/2)$ and a
$\zeta\in C^{\infty}_{0}(\bR^{d})$ with support in the ball of
radius $2\varepsilon$ centered at the origin and such that
$\zeta(x)=1$ for $|x|\leq\varepsilon$. Also take  a point
$(t_{0},x_{0})\in\bR^{d}_{T}$ and introduce $x(t)$ as a  solution (not
necessarily unique) of the problem
$$
x(t)=x_{0}+\int_{t_{0}}^{t}b(s,x(s))\,ds,\quad t \in\bR,
$$
where $b(t,x)$ is the vector with coordinates $b^{i}(t,x),i=1,
\ldots, d $.   Such solution exists since $b$ is locally bounded,
continuous in $x$ uniformly in $t$, and grows at most linearly in
$x$.

Set
$$
a_{0}^{ij}(t)=a^{ij}(t,x(t)),\quad b_{0}(t)=b(t,x(t)) ,\quad
c_{0}(t)=c(t,x(t)),
$$
$$
L_{0}=a_{0}^{ij}(t)D_{ij}+b_{0}(t)D_{i}-c_{0}(t), $$ $$ \quad
f_{0}(t)=f(t,x(t)),
$$
$$
u_{0}(t)=-\int_{t}^{\infty}f_{0}(s) \exp(-\int_{t}^{s }c_{0}(
r)\,dr)\,ds,
$$
$$
\eta(t,x)=\zeta(x-x(t)), \quad v(t,x)=[u(t,x) -u_{0}(t)]\eta(t,x)
.
$$

Observe that if   $\eta(t,x)\ne0$, then
$|x-x(t)|\leq2\varepsilon$, so that

$$
|a^{ij}(t,x)-a_{0}^{ij}(t)|\leq2^{\alpha}K\varepsilon^{\alpha},
\quad |b(t,x)-b_{0}(t)|\leq2^{\alpha}K\varepsilon^{\alpha}d,
$$
\begin{equation}
                                            \label{1.20.1}
|c(t,x)-c_{0}(t)|\leq2^{\alpha}K\varepsilon^{\alpha}, \quad
|f(t,x)-f_{0}(t)|\leq2^{\alpha}F_{\alpha}\varepsilon^{\alpha}.
\end{equation}
Also
$$
\eta_{t}(t,x)+b^{i}_{0}(t)\eta_{x_{i}}(t,x)=0, $$ $$ \quad
u_{0t}+L_{0}u_{0}= f_{0}.
$$

Next, by Lemma \ref{lemma 1.10.1} applied to
 $v$,  for   $x \in \bR^d $ such that
 $
|x-x_{0}|\leq\varepsilon 
 $,
we have $\eta(t_{0},x)=1$ and
$$
I:=\frac{|D^2u(t_{0},x)- D^2u (t_{0},x_{0})|}{|x-x_{0}|^{\alpha}}=
\frac{|D^2v(t_{0},x)- D^2v (t_{0},x_{0})|}{|x-x_{0}|^{\alpha}}
$$
$$
 \leq
N\sup_{s>t_{0}} [(v_{t}+L_{0}v)(s,\cdot)]_{\alpha}.
$$

Here
$$
v_{t}+L_{0}v =
\eta(u_{t}+L_{0}u-f_{0})+(u-u_{0})(\eta_{t}+L_{0}\eta+c_{0}\eta)
+2a^{ij}_{0}\eta_{x^{i}}u_{x^{j}}
$$
$$
=\eta (f-f_{0})+ \eta (L_{0}-L
)u+(u-u_{0})a^{ij}_{0}\eta_{x^{i}x^{j}}
+2a^{ij}_{0}\eta_{x^{i}}u_{x^{j}} .
$$

Since on the support of $\eta$ we have \eqref{1.20.1}, it is
standard to see that
$$
I\leq N(\varepsilon)F_{\alpha}+N\varepsilon^{\alpha}\sup_{s>t_{0}}
[u(s,\cdot)]_{2+\alpha}
$$
$$
+N(\varepsilon)\sup_{s>t_{0}}\|u(s,\cdot)\|_{2}
+N(\varepsilon)\sup|u_{0}|,
$$
where $N=N(d,\alpha,K,\delta)$ and  $N(\varepsilon)=N(
\varepsilon,d,\alpha,K,\delta)$. Due to the arbitrariness of
$x_{0}$ and $x$  and the fact that  obviously
   $|u_{0}(t)| \leq F_{0}$   in $\bR$,  we
obtain
$$
[u(t_{0},\cdot)]_{2+\alpha}\leq N(\varepsilon)(F_{0}+F_{\alpha})
$$
$$
+N\varepsilon^{\alpha}\sup_{s>t_{0}} [u(s,\cdot)]_{2+\alpha}
+N(\varepsilon)\sup_{s>t_{0}}\|u(s,\cdot)\|_{2}.
$$
Upon taking the sup of both sides with respect to $t_{0}\ge  T$ we
conclude
$$
\sup_{t \ge  T}[u(t,\cdot)]_{2+\alpha}\leq
N(\varepsilon)(F_{0}+F_{\alpha})
$$
$$
+N\varepsilon^{\alpha}\sup_{t\ge  T} [u(t,\cdot)]_{2+\alpha}
+N(\varepsilon)\sup_{t\ge  T}\|u(t,\cdot)\|_{2}.
$$
After choosing $\varepsilon$ appropriately, we finally get
\eqref{2.7.1}. The theorem is proved.

\begin{corollary}
                                            \label{corollary 4.15.3}
Let $T_{1}$, $T_{2}$, and $R $ be finite numbers such that $T\leq
T_{1}<T_{2}$ and    $R>0$. 
Then under
   Hypotheses \ref{hy} and \ref{H2} 
for any   $u\in\cC^{2+\alpha}(T)$ satisfying \eqref{2.6.4}
we have:  
$$|t-s|^{-1}
\, |u(s,x)-u(t,x)|+  |t-s|^{- \frac{1+ \alpha}{2}} |Du(s,x)- D
u(t,x)|
$$
\begin{equation}                                      \label{tt}
+ \, |t-s|^{- \alpha/2} \, |D^2u(s,x)- D^2u (t,x)|\leq N,
\end{equation}
  whenever $|x|\leq R$ and $s,t\in[T_{1},T_{2}]$,  
$t<s$,  where $N$ depends only on $\delta$, $\alpha$, $d$, $K$,
$F_{0}$, $F_{\alpha}$, $R$, and sup norms of
$|b(t,0)|,c(t,0),|f(t,0)|$ over $[T_{1},T_{2}]$.
\end{corollary}
    Proof.
   Take   $\zeta\in C^{\infty}_{0}(\bR^{d+1})$ such that
$\zeta(t,x)=1$ if  $t\in[T_{1},T_{2}]$ and $|x|\leq R$ and observe
that for $v:=u\zeta$  we have
$$
v_{t}=u\zeta_{t}-\zeta Lu+\zeta f.
$$
By combining this with \eqref{4.15.3} we see that
for any finite $t\geq T$
$$
 \|v_{t}(t,\cdot)\|_{\alpha}+ \|v(t,\cdot) \|_{2+\alpha}\leq N,
$$
where $N$ is as in the statement of the corollary. After that we
get our assertion from  Corollary \ref{corollary 4.13.2} and
Remark \ref{remark do}.

The a priori estimate \eqref{4.15.3} and the solvability of the
heat equation (Lemma \ref{lemma 4.14.1}) allow us to prove the
following result by the method of continuity.

\begin{lemma}
                                                   \label{lemma 4.15.3}
Under Hypotheses \ref{hy} and \ref{H2} additionally assume that
$b$, $c$, and $f$ are bounded. Then there exists a unique
$u\in\cC^{2+\alpha}(T)$ satisfying \eqref{2.6.4}.
\end{lemma}

Proof. Uniqueness follows from Theorem \ref{theorem 1.10.1}. In
the proof of existence first assume that $T>-\infty$ and take a
finite $S>T$. For $S$ and $T$ fixed introduce the space
$\frC^{2+\alpha} $ as a subspace of $\cC^{2+\alpha}(T)$ such that
$u(t,x)=0$ for $t\geq S$ and
$$
\|u\|_{\frC^{2+\alpha} }:=\sup_{t\in[T,S]
}\|u_{t}(t,\cdot)\|_{\alpha} +\sup_{t\in[T,S] }\|u
(t,\cdot)\|_{2+\alpha}<\infty.
$$
Also let $\frC^{\alpha}$ be the set of measurable functions
$g(t,x)$ such that $g(t,x)=0$ for $t\geq S$ and
$$
\|g\|_{\frC^{ \alpha} }:=\sup_{t\in[T,S] }\|g(t,\cdot)\|_{\alpha}
<\infty.
$$
One can easily check that $\frC^{2+\alpha}$ and $\frC^{\alpha}$
are Banach spaces.

For $\lambda\in[0,1]$ consider the following family of equations
\begin{equation}
                                                       \label{4.14.61}
u_{t}+[\lambda L+(1-\lambda)(\Delta -\delta)] u=g
\end{equation}
in $\bR^{d+1}_{T}$. We call a $\lambda\in[0,1]$ ``good" if for any
$g\in \frC^{\alpha}$ there is a unique solution $u
\in\frC^{2+\alpha}$ of \eqref{4.14.61}. Notice that if
$u\in\cC^{2+\alpha}(T)$ satisfies \eqref{4.14.61} with a
$g\in\frC^{\alpha}$, then by Theorem \ref{theorem 1.10.1} (with
$S$ in place of $T$) we have $u(t,x)=0$ for $t\geq S$ and
$$
   \sup_{t\ge  T}\|u(t,\cdot)\|_{2+\alpha} \leq N
\|g\|_{\frC^{\alpha} }.
$$
  From equation \eqref{4.14.61} we get an estimate of $u_{t}$ and
conclude that
\begin{equation}
                                                       \label{4.14.9}
   \|u\|_{\frC^{2+\alpha}}\leq N  \|g\|_{\frC^{\alpha}},
\end{equation}
where $N$ is independent of $\lambda$, $g$, $T$, and $S$.
Furthermore, by Lemma \ref{lemma 4.14.1} we have that $0$ is a
``good" point.

We now claim that all points of $[0,1]$ are ``good". To prove the
claim we take a ``good" point $\lambda_{0}$ (say $\lambda_{0}=0$)
and rewrite \eqref{4.14.61} as
\begin{equation}
                                                       \label{4.14.7}
u_{t}+[\lambda_{0} L+(1-\lambda_{0})\Delta] u=g
+(\lambda-\lambda_{0})(\Delta-L)u.
\end{equation}
Now fix $g\in \frC^{\alpha}$ and define a mapping $\cR$ which sends
$v\in
\frC^{2+\alpha}$ into the solution $u\in \frC^{2+\alpha}$ of the
equation
\begin{equation}
                                                       \label{4.14.8}
u_{t}+[\lambda_{0} L+(1-\lambda_{0})\Delta] u=g
+(\lambda-\lambda_{0})(\Delta-L)v.
\end{equation}
Observe that owing to our assumptions and the choice of
$\lambda_{0}$ the right-hand side of \eqref{4.14.8} is in
$\frC^{\alpha}$ and the mapping $\cR$ is well defined.

Estimate \eqref{4.14.9} shows that for any $v,w\in\frC^{2+\alpha}$
$$
\|\cR v-\cR w\|_{\frC^{2+\alpha}}\leq N|\lambda-\lambda_{0}| \| v-
w\|_{\frC^{2+\alpha}}
$$
with $N$ independent of $\lambda_{0},v$, and $w$. It follows that
there is an $\varepsilon>0$ such that for
$|\lambda-\lambda_{0}|\leq\varepsilon$ the mapping $\cR$ is a
contraction in $\frC^{2+\alpha}$ and has a fixed point $u$ which
obviously satisfies \eqref{4.14.7} and \eqref{4.14.61}. Therefore
such $\lambda$'s are ``good", which certainly proves our claim.

Owing to the boundedness of $f$, we have $fI_{[T,S]}\in
\frC^{\alpha}$, so that we now know that \eqref{2.6.4} has a
unique solution $u^{S} \in \cC^{2+\alpha}(T)$ if we take
$f(t,x)I_{t\le S}$ in place of $f$. By the above
\begin{equation}
                                                       \label{4.14.91}
\sup_{t\geq T}\|u^{S}(t, \cdot )\|_{2+\alpha}  \leq N\sup_{t\geq
T}\|f(t,\cdot)\|_{\alpha},
\end{equation}
   with $N$ independent of $f$, $T$, and $S$. This estimate  and
Corollary \ref{corollary 4.15.3} show  that the family
$u^{S},u^{S}_{x},u^{S}_{xx}$ is uniformly bounded and
equicontinuous on any bounded subset of $\bR^{d+1}_{T}$. By the
Ascoli-Arzel\`a theorem there is a continuous function $u$ on
$\bR^{d+1}_{T}$ having bounded and continuous derivatives 
with respect to $x$
up to
the second order and a sequence $S_{n}\to\infty$ such that
   $$
   (u^{S_{n}},u^{S_{n}}_{x},u^{S_{n}}_{xx}) \to(
u^{},u^{}_{x},u^{}_{xx})
$$
as $n \to \infty$,  uniformly on  bounded subsets of
$\bR^{d+1}_{T}$. Passing to the limit in the equations (of course,
in the integral form, see \eqref{2.7.5}) corresponding to
   $u^{S_n}$
shows that $u$ satisfies \eqref{2.6.4}. Furthermore,
\eqref{4.14.91} implies that with the same $N$
$$
\sup_{t\geq T}\|u(t,\cdot)\|_{2+\alpha} \leq N\sup_{t\geq
T}\|f(t,\cdot)\|_{\alpha}.
$$
Hence, $u\in \cC^{2+\alpha}(T)$ and we have proved the lemma if
$T>-\infty$.

To consider the case that $T=-\infty$ we define $u_{n}$ as
solutions of class $\cC^{2+\alpha}(-n)$ of \eqref{2.6.4} in
$\bR^{d+1}_{-n}$. Since the right-hand sides of the equations for
$u_{m}$ agree on $\bR^{d+1}_{-n}$ for $m\geq n$ by Theorem
\ref{theorem 1.10.1} we have that $u_{m}=u_{n}$ on
$\bR^{d+1}_{-n}$. Therefore, the limit $u$ of $u_{n}$ as $n\to\infty$
exists and obeys estimate \eqref{4.15.3} for $t\geq-n$ 
with arbitrary  $n$ and thus for any $t$. This
  shows that
$u\in\cC^{2+\alpha}(T)$ and finishes the proof of the lemma.

Here is a classical  result from the standard $C^{2+\alpha}$
theory of parabolic equations.

\begin{corollary}
                                                \label{corollary 4.15.1}
 In addition to the assumptions of Lemma \ref{lemma 4.15.3} suppose
 that $a,b,c$, and $f$ are $\alpha/2$-H\"older continuous in $t$
 with constant independent of $x$. Let $u \in \cC^{2+\alpha}(T)$
  be the solution
   of \eqref{2.6.4}.  Then,

(i)  For any $x \in \bR^d$, finite $ s, t$,
   with $T \le s < t $, estimate \eqref{tt} holds with a
   constant $N$ depending only  on
 $\delta$, $\alpha$, $d$,  $K$, $F_{\alpha}$, and
   the sup norms of $b$, $c$, and
 $f$;

(ii) The function
   $u_{t}$ is bounded on $\bR^{d+1}_{T}$,
 $\alpha$-H\"older continuous in $x$ with constant independent of
 $t$ and $\alpha/2$-H\"older continuous in $t$ with constant
independent of $x$.
\end{corollary}
 
 Indeed, in    assertion (ii) the
  the
 $\alpha $-H\"older continuity in $x$
of
$u_{t}$   is obtained
directly from the equation  even without requiring the H\"older
continuity of the data in $t$.

The remaining assertions of the lemma
  follow  directly from the equation due to 
Theorem \ref{theorem 1.10.1},
 Corollary
\ref{corollary 4.13.2},
  and   Remark \ref{remark do}.

One could easily provide   estimates of the above mentioned
H\"older constants. We leave this to the reader.

{\bf Proof of Theorem \ref{theorem 2.7.1}}.  As usual uniqueness
follows from the a priori estimate (see  Theorem \ref{theorem
1.10.1}). To prove existence we use truncations. Introduce
$\chi_{n}(\tau)=(-n)\vee \tau\wedge n$, $n = 1,2,...$,
   $b^{i}_{n}=\chi_{n}(b^{i})$, $i = 1, \ldots, d$, and similarly
introduce $c_{n}$ and $f_{n}$. Since $\chi_{n}$ are Lipschitz
continuous with constant one and, for any $\tau,t\geq0$, we have
$\chi_{n}(\tau t)\leq\tau\chi_{n}(t)$, the truncated coefficients
satisfy Hypothesis \ref{hy} with the same constants $K,F_{0},$ and
$F_{\alpha}$.

By Lemma \ref{lemma 4.15.3} there exists  a unique
$u_{n}\in\cC^{2+\alpha}(T)$ satisfying equation \eqref{2.6.4} with
$b_{n},c_{n}$, and $f_{n}$ in place of $b,c$, and $f$,
respectively. By Theorem \ref{theorem 1.10.1}
\begin{equation}
                                            \label{2.7.6}
   \sup_{t\ge  T}\|u_{n}(t,\cdot)\|_{2+\alpha} \leq N (F_0
+ F_{\alpha}),
\end{equation}
where $N$ is independent of $n$. 
This estimate  and Corollary \ref{corollary 4.15.3} 
imply that the family
$u_{n},u_{nx},u_{nxx}$ is uniformly bounded and equicontinuous in any
bounded subset of $\bR^{d+1}_{T}$.

Since, for any bounded set $\Gamma\subset  \bR^{d+1}_{T}$ there
exists an $m$ such that $(b_{n},c_{n},f_{n})=(b,c,f)$ on $\Gamma$
for all $n\geq m$, to prove the theorem it only remains to repeat
the part of the proof of Lemma \ref{lemma 4.15.3} concerning
$u^{S_{n}}$. The theorem is proved.

\end{document}